\documentstyle[amsfonts,12pt]{article}
\input amssym.def

\newcommand{\bul}{{\bullet}}
\newcommand{\al}{{\alpha}}

\newcommand{\mg}{{\mathfrak{g}}}
\newcommand{\mh}{{\mathfrak{h}}}
\newcommand{\mv}{{\mathfrak{v}}}
\newcommand{\mP}{{\mathfrak{P}}}

\newcommand{\cU}{U^f}

\newcommand{\mmU}{{\mathfrak{U}}}

\newcommand{\mU}{{\cal U}}

\newcommand{\bmU}{{\underline{\cal U}}}

\newcommand{\tmU}{{\tilde{\cal U}}}

\newcommand{\bU}{{\bf U}}

\newcommand{\bQ}{{\bf Q}}

\newcommand{\Om}{{\Omega}}

\newcommand{\si}{{\sigma}}
\newcommand{\ga}{{\gamma}}

\newcommand{\ve}{{\varepsilon}}

\newcommand{\G}{{\Gamma}}
\newcommand{\pa}{{\partial}}

\newcommand{\cF}{{\cal F}}
\newcommand{\cD}{{\cal D}}
\newcommand{\T}{{\cal T}}
\newcommand{\cZ}{{\cal Z}}

\newcommand{\cB}{{\cal B}}

\newcommand{\bbR}{{\Bbb R}}
\newcommand{\bbZ}{{\Bbb Z}}
\newcommand{\bbRf}{{{\Bbb R}^d_{formal}}}

\newcommand{\n}{{\nabla}}

\newcommand{\de}{{\delta}}
\newcommand{\D}{{\Delta}}

\newcommand{\tF}{{\tilde{F}}}
\newcommand{\SM}{{\cal S}M}
\newcommand{\OmS}{\Om(M,\SM)}
\newcommand{\OmT}{\Om(M,\T_{poly})}
\newcommand{\OmD}{\Om(M,\cD_{poly})}
\newcommand{\FT}{\cF^0\T_{poly}}
\newcommand{\FD}{\cF^0\cD_{poly}}

\date{}
\newtheorem{defi}{Definition}
\newtheorem{pred}{Proposition}

\newtheorem{teo}{Theorem}
\newtheorem{cor}{Corollary}
\begin{document}

\begin{center}
{\Large\bf Covariant and Equivariant Formality Theorems.}\\[0.5cm]
Vasiliy Dolgushev\footnote{On leave of absence from: Univ. Center of JINR (Dubna) and
ITEP (Moscow)} \\[0.3cm]
{\it Department of Mathematics, MIT,} \\
{\it 77 Massachusetts Avenue,} \\
{\it Cambridge, MA, USA 02139-4307,}\\
{\it E-mail address: vald@math.mit.edu}
\end{center}

\begin{abstract}
We give a proof of Kontsevich's formality theorem
for a general ma\-ni\-fold using Fedosov resolutions of
algebras of polydifferential operators and polyvector
fields. The main advantage of our construction
of the formality quasi-isomorphism is that it is based on the
use of covariant tensors unlike Kontsevich's original proof,
which is based on $\infty$-jets
of polydifferential operators and polyvector fields.
Using our
construction we prove that if a group $G$ acts smoothly on a
ma\-ni\-fold $M$ and $M$ admits a $G$-invariant affine connection
then there exists a $G$-equivariant quasi-isomorphism of formality.
This result implies that if a ma\-ni\-fold $M$ is equipped with
a smooth action of a finite or compact group $G$ or
equipped with a free action of a Lie group $G$ then
$M$ admits a $G$-equivariant formality quasi-isomorphism.
In particular, this gives a solution of the
deformation quantization problem for
an arbitrary Poisson orbifold.\\[0.3cm]
MSC-class: 16E45; 53C15; 18G55.
\end{abstract}

\section{Introduction}
Preserving symmetries in quantization procedures is
one of the most important problems in mathematical physics.
In this paper we discuss the problem
for Kontsevich's  formality quasi-iso\-mor\-phism
 \cite{Kontsevich}.

The purpose of this paper is twofold. First, we propose
a manifestly covariant construction of Kontsevich's
formality quasi-iso\-mor\-phism for a general
smooth ma\-ni\-fold $M$ using the Fedosov resolutions of the
algebras of po\-ly\-dif\-fe\-ren\-tial operators and po\-ly\-vec\-tor fields
and Kontsevich's quasi-iso\-mor\-phism of formality for
the space $\bbR^d$.
The quasi-iso\-mor\-phism obtained by our procedure depends on
an affine torsion free connection on $M$.
Second, we consider a ma\-ni\-fold $M$ with a smooth action
of a group $G$ and show that if $M$ admits
a $G$-invariant affine connection then the
formality quasi-iso\-mor\-phism corresponding to
the connection is $G$-equivariant.

Our method may be regarded as a generalization
of the work \cite{CFT}, in which a covariant construction
of a star-product on an arbitrary Poisson ma\-ni\-fold is presented.
As in \cite{CFT}, we use the Fedosov
resolution\footnote{In work \cite{CF} construction of the resolution
is called Weinstein's exponential map \cite{Wein}.}
of the algebra of functions
on $M$ in the algebra of sections of
the formally completed symmetric algebra of the
cotangent bundle $T^*M$ and a fiberwise formality
quasi-iso\-mor\-phism. However,
we use a fixed Fedosov differential and instead modify
the fiberwise formality quasi-iso\-mor\-phism
unlike the authors of \cite{CFT}.
It is this modification that allows us to get a more
general result, namely, to construct a formality
quasi-iso\-mor\-phism for an arbitrary
smooth ma\-ni\-fold.

A sketchy proof for the formality of the algebra of
po\-ly\-dif\-fe\-ren\-tial operators on an arbitrary smooth
ma\-ni\-fold is originally given in \cite{Kontsevich} (see section
7). A more detailed explanation of the proof can be found in the
appendix of paper \cite{K}. In this context it is also worth
mentioning paper \cite{BLN}, in which a covariant deformation quantization
for an arbitrary Poisson manifold is proposed
via a path integral approach.

Our construction of the formality quasi-iso\-mor\-phism
for a general ma\-ni\-fold $M$ may be regarded as a modification of
Kontsevich's approach \cite{Kontsevich}, \cite{K} in which
the $\infty$-jets are replaced by infinite collections of
symmetric covariant tensors and the flat connection
on the jet bundle is replaced by Fedosov differential,
constructed with the help of some torsion free
connection on $M$. This modification turns out to be much simpler
for applications because at all the stages of
our construction of the quasi-iso\-mor\-phism we deal with manifestly
covariant objects.

We use our construction of the formality quasi-iso\-mor\-phism
in order to prove an equivariant formality theorem (see theorem 5 in
section 5), which is the main result of the paper. This theorem allows us
to get several interesting corollaries. Namely, it turns
out that if a finite or compact group $G$ acts smoothly
on a ma\-ni\-fold $M$ or a Lie group $G$ acts freely on $M$ then
$M$ admits a $G$-equivariant formality
quasi-iso\-mor\-phism. In particular, this gives us a solution of
the deformation quantization problem for an arbitrary
Poisson orbifold (see corollary 3 in section 5).

The structure of this paper is as follows.
The next section is devoted to basic notations we use throughout
the paper. In this section we recall some required notions of
homotopy theory of differential graded Lie algebras, we review
necessary properties of Kontsevich's formality quasi-iso\-mor\-phism
\cite{Kontsevich} for the space $\bbR^d$, and finally we introduce
some notions required for a construction of the Fedosov resolutions
of algebras of po\-ly\-vec\-tor fields and
po\-ly\-dif\-fe\-ren\-tial operators.
In the third section we present the Fedosov
resolutions of the algebras of functions, po\-ly\-dif\-fe\-ren\-tial
operators, and po\-ly\-vec\-tor fields on a smooth ma\-ni\-fold.
Section 4 is devoted to our construction of
a formality quasi-iso\-mor\-phism for a general smooth
ma\-ni\-fold. In section 5 we prove an equivariant formality
theorem and present its corollaries. Finally, in the concluding
section of the paper we discuss possible applications
and generalizations of the equivariant formality theorem.

Throughout the paper the summation over repeated indices
is assumed. Sometimes we omit the prefix ``super-'' referring to
super-algebras, Lie super-brackets, and su\-per(co)\-com\-mu\-ta\-tive
(co)multiplications. We assume that $M$ is a smooth real ma\-ni\-fold
of dimension $d$\,. We omit symbol $\wedge$ referring to a local basis
of exterior forms as if we thought of $dx^i$'s as anti-commuting variables.
Finally, we always assume that a nilpotent linear operator is
the one whose second power is vanishing.

\section{Preliminaries.}
This section is devoted to basic notations we use throughout
the paper. Although the objects we introduce here are well known
we assemble in this section all basic definitions and results we need since
notations and terminology vary from source to another.

We start with a sketchy introduction of some required notions
of homotopy theory of differential graded Lie algebras (DGLA).
A more detailed discussion of the theory can be found in
paper \cite{HS}.

Let $(\mh, d_{\mh}, \{,\})$ be a DGLA. We assume that the
differential $d_{\mh}$ is of degree one and the Lie super-bracket
$\{,\}$ is of degree zero. To include $(\mh, d_{\mh}, \{,\})$ in the
context of strong homotopy Lie algebras we associate to
 $\mh$ a coassociative cocommutative
coalgebra $C_{\bul}(\mh[1])$ cofreely cogenerated by the vector space $\mh$
with a shifted parity. The DGLA structure $(d_{\mh}, \{,\})$
induces on $C_{\bul}(\mh[1])$ a coderivation $Q$ with
two non-vanishing structure maps
\begin{equation}
\label{Q-str}
\begin{array}{c}
Q_1=d_{\mh}\,:\, \mh \mapsto \mh[1]\,, \\[0.3cm]
Q_2=\{ , \} \,:\, \wedge^2 \mh \mapsto \mh\,.
\end{array}
\end{equation}
With this definition the set of axioms of the DGLA $(\mh, d_{\mh}, \{,\})$ is
equivalent to nilpotency of coderivation (\ref{Q-str})
$$
Q^2=0\,.
$$

Given two DG Lie algebras $(\mh_1, d_1, [,]_1)$ and $(\mh_2, d_2, [,]_2)$
we shall be interested in their morphisms in the category of
strong homotopy Lie algebras, namely

\begin{defi}
An $L_{\infty}$-morphism $F$ from the DGLA
$(\mh_1, d_1, [,]_1)$ to the DGLA $(\mh_2, d_2, [,]_2)$ is
a homomorphism of the
coassociative cocommutative coalgebras
\begin{equation}
\label{U=homo}
F\,:\,C_{\bul}(\mh_1[1]) \mapsto C_{\bul}(\mh_2[1])
\end{equation}
compatible with the nilpotent coderivations $Q_1$ and
$Q_2$ corresponding to the DGLA structures $(d_1, [,]_1)$
and  $(d_2, [,]_2)$, respectively
$$
Q_2 F (X)= F(Q_1 X)\,, \qquad \forall~X\in C_{\bul}(\mh_1[1])\,.
$$
\end{defi}

Furthermore,

\begin{defi}
A quasi-iso\-mor\-phism $F$ from the DGLA
$(\mh_1, d_1, [,]_1)$ to the DGLA $(\mh_2, d_2, [,]_2)$
is an $L_{\infty}$-morphism from
$\mh_1$ to $\mh_2$ whose first
structure map
$$
F_1 : \mh_1 \mapsto \mh_2
$$
gives an iso\-mor\-phism of the spaces of cohomologies
$H^{\bul}(\mh_1,d_1)$ and $H^{\bul}(\mh_2,d_2)$.
\end{defi}
In what follows the notation
$$
F\,:\, (\mh_1, d_1, [,]_1) \leadsto (\mh_2, d_2, [,]_2)
$$
means that $F$ is a quasi-iso\-mor\-phism form the
DGLA $(\mh_1, d_1, [,]_1)$ to the DGLA
$(\mh_2, d_2, [,]_2)$\,.

By unfolding the formal definition we see that
a homomorphism $F$ of coassociative cocommutative coalgebras
$C_{\bul}(\mh_1[1])$ and  $C_{\bul}(\mh_2[1])$
is uniquely defined by a semi-infinite collection
of polylinear maps
\begin{equation}
\label{struct}
F_n : \wedge^n \mh_1 \mapsto \mh_2 [1-n], \qquad n\ge 1
\end{equation}
and the compatibility of $F$ with the coderivations $Q_1$ and
$Q_2$ on $C_{\bullet}(\mh_1[1])$ and $C_{\bullet}(\mh_2[1])$,
respectively, is equivalent to the following semi-infinite
collection of equations

$$
d_2 F_n(\ga_1, \ga_2, \ldots, \ga_n)-
\sum_{i=1}^n (-)^{k_1+\ldots+k_{i-1}+1-n}
F_n(\ga_1, \ldots, d_1 \ga_i, \ldots, \ga_n)=
$$
\begin{equation}
=\frac12 \sum_{k,l\ge 1,~ k+l=n} \frac1{k!l!}
 \sum_{\ve\in S_n}
\pm [F_k (\ga_{\ve_1}, \ldots, \ga_{\ve_k}), F_l (\ga_{\ve_{k+1}}, \ldots,
\ga_{\ve_{k+l}})]_2-
\label{q-iso}
\end{equation}
$$
-\sum_{i\neq j}
\pm F_{n-1}([\ga_i,\ga_j]_{1}, \ga_1, \ldots, \hat{\ga_i}, \ldots, \hat{\ga_j}, \ldots \ga_n),
\qquad \ga_i \in \mh_1^{k_i}\,.
$$
\vspace{0.1cm}
{\bf Remark.} Notice that in order to define the signs in formulas (\ref{q-iso})
one should use a rather complicated rule.
For example, the signs that stand before
the terms of the first sum at the right hand side depend
on permutations $\ve\in S_n$, on
degrees of $\ga_i$, and on the numbers $k$ and $l$.
The simplest way to check that all the signs are correct
is to show that the right hand side of equation (\ref{q-iso})
is closed with respect to the following
differential acting on the space of graded polylinear maps
$$
d_{Hom} \,:\, Hom (\wedge^n \mh_1, \mh_2[k]) \mapsto Hom (\wedge^n \mh_1, \mh_2[k+1]),
$$
\begin{equation}
d_{Hom} \Psi(\ga_1, \ga_2, \ldots, \ga_n)= d_2 \Psi(\ga_1, \ga_2, \ldots, \ga_n)-
\label{d-Hom}
\end{equation}
$$
-\sum_{i=1}^n (-)^{k_1+\ldots+k_{i-1}+ k}
\Psi(\ga_1, \ldots, d_1 \ga_i, \ldots, \ga_n), \qquad \ga_i \in \mh^{k_i}_1,
$$
where $\Psi\in Hom (\wedge^n \mh_1, \mh_2[k])$\,.

~\\
{\bf Example.} An important example of a quasi-iso\-mor\-phism
from a DGLA $\mh_1$ to a DGLA $\mh_2$ is provided by a
DGLA-homomorphism
$$
H\,:\,\mh_1\mapsto \mh_2\,,
$$
which induces an iso\-mor\-phism on the spaces
of cohomologies
$H^{\bul}(\mh_1,d_1)$ and $H^{\bul}(\mh_2,d_2)$.
In this case the quasi-iso\-mor\-phism has the only
non-vanishing structure map
$$
F_1=H\,.
$$

To conclude the introductory part on the homotopy theory
we recall that a DGLA algebra $(\mh,d, [,])$ is called 
formal\footnote{Our definition of formality is slightly 
different from the conventional one. However, if $\mh$ is formal in the 
sense of our definition then it is conventionally formal.} 
if there is a quasi-isomorphism from the 
graded Lie algebra $H^{\bul}(\mh,d)$ to the DGLA $\mh$.

Let now $M$ be a smooth ma\-ni\-fold of dimension $d$ and
$D_{poly}(M)$ be a vector space of po\-ly\-dif\-fe\-ren\-tial operators
on $M$

\begin{equation}
\label{D-poly}
D_{poly}(M)=\bigoplus_{k=-1}^{\infty} D_{poly}^k(M) \,,
\end{equation}
where the $D_{poly}^k(M)$ consists of po\-ly\-dif\-fe\-ren\-tial
operators of rank $k+1$
$$
\Phi\,:\,C^{\infty}(M)^{\otimes (k+1)} \mapsto C^{\infty}(M)\,.
$$

The space $D_{poly}(M)$ can be endowed with the so-called
Gerstenhaber bracket which is defined
between homogeneous elements $\Phi_1\in D^{k_1}_{poly}(M)$ and
$\Phi_2\in D^{k_2}_{poly}(M)$ as follows

$$[\Phi_1, \Phi_2](a_0,\,\ldots, a_{k_1+k_2})=$$
\begin{equation}
\label{Gerst}
\sum_{i=0}^{k_1}(-)^{ik_2}
\Phi_1(a_0,\,\ldots, a_{i-1},\Phi_2 (a_i,\,\ldots,a_{i+k_2}),\, \ldots, a_{k_1+k_2})
\end{equation}
$$
-(-)^{k_1k_2} (1 \leftrightarrow 2)\,.
$$

Direct computation shows that (\ref{Gerst}) is a Lie (super)bracket and
therefore $D_{poly}$ is Lie (super)algebra.
The multiplication operator $m_0\in D_{poly}^1(M)$
$$m_0 : C^{\infty}(M) \otimes C^{\infty}(M)\mapsto C^{\infty}(M)$$
satisfies the associativity condition, which can be written
in terms of bracket (\ref{Gerst}) as

\begin{equation}
\label{assoc}
[m_0,m_0]=0\,.
\end{equation}
Thus $m_0$ defines a nilpotent interior derivation of $D_{poly}(M)$

\begin{equation}
\label{pa}
\pa \Phi = [m_0,\Phi] : D^k_{poly}(M)\mapsto D^{k+1}_{poly}(M)\,, \qquad
\pa^2 =0\,,
\end{equation}
which turns the Lie algebra $D_{poly}(M)$ into a DGLA.

A DGLA of cohomologies of $D_{poly}(M)$ is described by the
Hochschild-Kostant-Rosenberg theorem, which says that
$$
H^{\bullet}(D_{poly}(M), \pa)= T_{poly}(M)\,,
$$
where $T_{poly}(M)$ is a DGLA of the po\-ly\-vec\-tor fields with a vanishing differential
$d_T=0\,:\, T^k_{poly}(M)\mapsto T^{k+1}_{poly}(M)$ and a Lie bracket being
the standard Schouten-Nijenhuis bracket. Namely,

\begin{equation}
\label{T}
T_{poly}(M)=\bigoplus_{k=-1}^{\infty} T^k_{poly}(M)\,, \qquad T^k_{poly}(M)=\G(\wedge^{k+1} TM)\,,
\end{equation}
and the Schouten-Nijenhuis bracket $[,]_{SN}$
is defined as an ordinary Lie bracket between vector fields
and then extended by Leibniz rule with respect to the $\wedge$-product
to an arbitrary pair of po\-ly\-vec\-tor
fields.

The formality theorem by Kontsevich \cite{Kontsevich} says that
for any smooth ma\-ni\-fold $M$ there exists a quasi-iso\-mor\-phism
from the DGLA $T_{poly}(M)$ of po\-ly\-vec\-tor fields on $M$ to the DGLA
$D_{poly}(M)$ of po\-ly\-dif\-fe\-ren\-tial operators on $M$. In our paper
we use this result for the case when $M=\bbR^d$\,.

In paper \cite{Kontsevich} Kontsevich proposed an interesting technique
for computing structure maps of a quasi-iso\-mor\-phism $U$ from the
DGLA $T_{poly}(\bbR^d)$ of po\-ly\-vec\-tor fields
to the DGLA $D_{poly}(\bbR^d)$ of po\-ly\-dif\-fe\-ren\-tial
operators on $\bbR^d$. Although the existence of the formality
quasi-iso\-mor\-phism for $\bbR^d$
has been also proved by Tamarkin \cite{Dima} here we
need explicit Kontsevich's construction
because the quasi-iso\-mor\-phism $U$ given
in \cite{Kontsevich} satisfy certain
peculiar properties\footnote{see the beginning of
section 7 in \cite{Kontsevich}}, which we use in
our construction of the formality quasi-iso\-mor\-phism
for a general ma\-ni\-fold.
We assemble the properties of $U$ in the following

\begin{teo}[Kontsevich, \cite{Kontsevich}] \label{1}
There exists a quasi-iso\-mor\-phism $U$

\begin{equation}
\label{form}
U\,:\,T_{poly}(\bbR^d)\leadsto D_{poly}(\bbR^d)
\end{equation}
from the
DGLA $T_{poly}(\bbR^d)$ of po\-ly\-vec\-tor fields
to the DGLA $D_{poly}(\bbR^d)$ of po\-ly\-dif\-fe\-ren\-tial
operators on the space $\bbR^d$ such that
\begin{enumerate}
\item One can replace $\bbR^d$ in (\ref{form}) by its formal completion $\bbRf$
at the origin.

\item The quasi-iso\-mor\-phism $U$ is equivariant with respect to
linear transformations of the coordinates on $\bbRf$\,.

\item If $n>1$ then

\begin{equation}
\label{vanish}
U_{n}(v_1, v_2, \ldots, v_n)=0
\end{equation}
for any set of vector fields
$v_1, v_2, \ldots, v_n\in T^0_{poly}(\bbRf)$\,.

\item If $n\ge 2$ and $v\in T^0_{poly}(\bbRf)$ is linear in the
coordinates on $\bbRf$ then for any set of po\-ly\-vec\-tor fields
$\ga_2, \ldots, \ga_n\in T_{poly}(\bbRf)$

\begin{equation}
\label{vanish1}
U_{n}(v,\ga_2, \ldots, \ga_n)=0\,.
\end{equation}
\end{enumerate}
\end{teo}

We now turn to some definitions required for
constructing Fedosov resolutions of algebras of po\-ly\-vec\-tor
fields and po\-ly\-dif\-fe\-ren\-tial operators.
First, we give a definition of a bundle $\SM$ of the formally
completed symmetric algebra of the cotangent bundle $T^*M$.
This bundle is a natural analogue of the Weyl algebra bundle
used in paper \cite{Fedosov} by Fedosov.

\begin{defi} The bundle $\SM$ of formally completed symmetric algebra
 of the cotangent bundle $T^*M$
is defined as a bundle over the ma\-ni\-fold $M$  whose sections
are infinite collections of symmetric covariant tensors $a_{i_1\ldots
i_p}(x)$\,, where $x^i$ are local coordinates, $p$
runs from $0$ to $\infty$\,, and the indices
$i_1,\ldots, i_p$ run from $1$ to $d$\,.
\end{defi}
It is convenient to introduce auxiliary variables $y^i$\,, which
transform as contravariant vectors. This allows us to
rewrite any section $a\in \G(\SM)$ in the form
of the formal power series

\begin{equation}
\label{sect}
a=a(x,y)=\sum_{p=0}^{\infty} a_{i_1\ldots
i_p}(x)y^{i_1}\ldots y^{i_p}\,.
\end{equation}
In this way the variables $y^i$ may be thought of as formal
coordinates on the fibers of the tangent bundle $TM$.

It is easy to observe that the vector space $\G(\SM)$ is
naturally endowed with the commutative product which is induced by a
fiberwise multiplication of formal power series in $y^i$\,. This product
makes $\G(\SM)$ into a commutative algebra with a unit.

Now we turn to definitions of
formal fiberwise po\-ly\-vec\-tor fields and formal fiberwise
po\-ly\-dif\-fe\-ren\-tial operators on $\SM$\,.

\begin{defi}
A bundle $\T^k_{poly}$ of formal fiberwise po\-ly\-vec\-tor
fields of degree $k$ is a bundle over $M$ whose sections
are $C^{\infty}(M)$-linear operators
$\mv : \wedge^{k+1} \G(\SM) \mapsto \G(\SM)$ of the form
\begin{equation}
\label{vect}
\mv =\sum_{p=0}^{\infty}\mv^{j_0\dots j_k}_{i_1\dots i_p}(x)y^{i_1}
\ldots y^{i_p} \frac{\pa}{\pa y^{j_0}}\wedge  \ldots \wedge \frac{\pa}{\pa
y^{j_k}}\,,
\end{equation}
where we assume that the infinite sum in $y$'s is formal and
$\mv^{j_0\dots j_k}_{i_1\dots i_p}(x)$ are tensors symmetric in
indices $i_1, \ldots, i_p$ and antisymmetric in indices
$j_0, \ldots, j_k$\,.
\end{defi}

Extending the definition of the formal fiberwise po\-ly\-vec\-tor
field by allowing the fields to be inhomogeneous
we define the total bundle $\T_{poly}$ of
formal fiberwise po\-ly\-vec\-tor fields

\begin{equation}
\label{cal-T}
\T_{poly} =\bigoplus_{k=-1}^{\infty} \T_{poly}^k\,, \qquad
\T_{poly}^{-1}=\SM\,.
\end{equation}

We mention that the fibers of the bundle $\T_{poly}$ form
a DGLA $T_{poly}(\bbRf)$ of po\-ly\-vec\-tor fields on
the formal completion $\bbRf$ of $\bbR^d$ at the origin.

\begin{defi}
A bundle $\cD^k_{poly}$ of formal fiberwise po\-ly\-dif\-fe\-ren\-tial
operator of degree $k$ is a bundle over $M$ whose sections
are $C^{\infty}(M)$-polylinear
maps $\mP : \bigotimes^{k+1} \G(\SM) \mapsto \G(\SM)$ of the form
\begin{equation}
\label{operr}
\mP =\sum_{\al_0 \ldots \al_k}\sum_{p=0}^{\infty}\mP^{\al_0\dots \al_k}_{i_1\dots i_p}(x)y^{i_1}
\ldots y^{i_p} \frac{\pa}{\pa y^{\al_0}}\otimes  \ldots \otimes \frac{\pa}{\pa
y^{\al_k}}\,,
\end{equation}
where $\al$'s are a multi-indices $\al={j_1\ldots j_l}$ and
$$
\frac{\pa}{\pa y^{\al}}=\frac{\pa}{\pa y^{j_1}} \ldots \frac{\pa}{\pa
y^{j_l}}\,,
$$
the infinite sum in $y$'s is formal, and the
sum in the orders of derivatives $\pa/\pa y$
is finite.
\end{defi}
Notice that the tensors
$\mP^{\al_0\dots \al_k}_{i_1\dots i_p}(x)$ are
symmetric in covariant indices $i_1,\ldots, i_p$\,.

As well as for po\-ly\-vec\-tor fields we define the total bundle
$\cD_{poly}$ of formal fiberwise po\-ly\-dif\-fe\-ren\-tial operators
as the direct sum

\begin{equation}
\label{cal-D}
\cD_{poly} =\bigoplus_{k=-1}^{\infty} \cD^k_{poly}\,, \qquad
\cD^{-1}_{poly}=\SM\,.
\end{equation}

We mention that the fibers of the bundle $\cD_{poly}$ form
a DGLA $D_{poly}(\bbRf)$ of po\-ly\-dif\-fe\-ren\-tial operators
on $\bbRf$\,.

For our purposes we need to tensor the bundles we introduced with the
exterior algebra bundle $\bigwedge T^*M$. Namely, instead of the commutative algebra
$\G(\SM)$ of formal power series in the fiber coordinates $y^i$ of
the tangent bundle $TM$ we will need a super-commutative algebra $\Om(M,\SM)$ of
exterior forms on $M$ with values in $\SM$, namely

\begin{equation}
\label{elem}
\Om(M,\SM)=
\{ a(x,y,dx)=\sum_{p,q\ge 0} a_{i_1\ldots
i_p j_1 \ldots j_q}(x)y^{i_1}\ldots y^{i_p}dx^{j_1} \ldots dx^{j_q} \}\,,
\end{equation}
where $a_{i_1\ldots i_p j_1 \ldots j_q}(x)$ are covariant tensors symmetric in indices
$i_1,\ldots, i_p$ and antisymmetric in indices $j_1,\ldots,
j_q$\,.

Algebra $\Om(M,\SM)$ is $\bbZ$-graded with respect to degree $p$ in ``$y$''
and $\bbZ$-graded with respect to the ordinary exterior
degree $q$.
$$
\Om(M,\SM)= \bigoplus_{p,q\ge 0} \Om^q(M,{\cal S}^p M)\,.
$$

Next, we introduce vector spaces $\OmT$
and $\OmD$ of smooth exterior forms
on $M$ with values in $\T_{poly}$ and $\cD_{poly}$ respectively.
It is easy to see that both $\OmT$ and $\OmD$ are naturally endowed with
DGLA structures induced by the respective fiberwise DGLA structures on
$T_{poly}(\bbRf)$ and $D_{poly}(\bbRf)$\,.
We denote the differential and the Lie bracket in $\OmD$
by $\pa$ and $[,]$\,, respectively, and denote the
Lie bracket in $\OmT$ by $[,]_{SN}$\,.

Although $\Om(M,\SM)$ is a commutative subalgebra
of the DGLA $\OmT$ ($\OmD$) of po\-ly\-vec\-tor fields
(po\-ly\-dif\-fe\-ren\-tial operators) of degree $-1$ we
consider $\Om(M,\SM)$ separately since we refer to $\Om(M,\SM)$
not only as to a Lie subalgebra of the DGLA $\OmT$ ($\OmD$)
but also as to a super-commutative algebra with
an ordinary multiplication.

The DG Lie algebras $\OmT$ and $\OmD$ are $\bbZ\times\bbZ\times\bbZ$-graded.
Namely, they have an obvious grading of the exterior forms, grading
with respect to the degree of po\-ly\-vec\-tor field or po\-ly\-dif\-fe\-ren\-tial
operator, and finally grading with respect to a difference
of degrees in $y$ and in $\pa/\pa y$\,.
While the parity of elements of $\Om(M,\SM)$ is defined by the exterior degree
$q$\,, the parity of the elements of $\OmT$ ($\OmD$) is defined by
the sum $r=q+k$ of the exterior degree $q$ and the degree $k$
of a po\-ly\-vec\-tor field (po\-ly\-dif\-fe\-ren\-tial operator).

Due to properties $1$ and $2$ in theorem \ref{1} we have
a fiberwise quasi-iso\-mor\-phism

\begin{equation}
\label{cal-U}
\cU :  (\OmT,0,[,]_{SN}) \leadsto (\OmD,\pa,[,])\,.
\end{equation}
from the DGLA $(\OmT,0,[,]_{SN})$ to the DGLA $(\OmD,\pa,[,])$\,.
We use the quasi-iso\-mor\-phism (\ref{cal-U}) in
section $4$ in order to prove the formality theorem for
a general smooth ma\-ni\-fold.
Now we turn to another important ingredient of our construction.

\section{Fedosov resolutions of $T_{poly}(M)$ and $D_{poly}(M)$}
In this section we construct resolutions for the DGLA $D_{poly}(M)$
of po\-ly\-dif\-fe\-ren\-tial operators and the DGLA $T_{poly}(M)$ of
po\-ly\-vec\-tor fields on an arbitrary smooth ma\-ni\-fold $M$ using the
DGLA $\OmD$ of exterior forms with values in fiberwise
po\-ly\-dif\-fe\-ren\-tial operators and the DGLA $\OmT$ of exterior forms
with values in fiberwise po\-ly\-vec\-tor fields on $\SM$. These
resolutions may be regarded as classical analogs of the
construction of the so-called quantum exponential Fedosov map from
the algebra of functions on a symplectic ma\-ni\-fold to an algebra of
sections of the Weyl bundle \cite{Fedosov}. In this sense the
resolutions are more reminiscent of what is called generalized
formal exponential map used in work \cite{CFT} and discussed in
papers \cite{Wein} and \cite{GRSh}.

We will proceed with the DG Lie algebras $D_{poly}(M)$ and $T_{poly}(M)$
and the algebra of functions $C^{\infty}(M)$ simultaneously and
denote the same operations on different algebras $\OmT$, $\OmD$
and $\Om(M,\SM)$  by the same letters. In what follows it does
not lead to any confusion.

The differential
\begin{equation}
\label{del}
\de= dx^i \frac{\pa}{\pa y^i} \,:\, \Om^q(M,\SM) \mapsto \Om^{q+1}(M,\SM)\,, \qquad \de^2=0
\end{equation}
on the algebra $\OmS$ obviously extends to differentials
on $\OmT$ and $\OmD$. Namely,

\begin{equation}
\label{delta}
\de= [dx^i \frac{\pa}{\pa y^i},\bullet\,]_{SN} \,:\, \Om^q (M, \T_{poly}) \mapsto
 \Om^{q+1} (M, \T_{poly})\,, \qquad \de^2=0\,,
\end{equation}
and

\begin{equation}
\label{delta1}
\de= [dx^i \frac{\pa}{\pa y^i},\bullet\,] \,:\,  \Om^q (M, \cD_{poly}) \mapsto
 \Om^{q+1} (M, \cD_{poly})\,, \qquad \de^2=0\,.
\end{equation}

Since the multiplication $m\in \G(\cD^1_{poly})$ in $\G(\SM)$
is $\de$-closed
$$
\de m =0
$$
$\de$ (anti)commutes with the differential $\pa$ in $\OmD$\,.
By definition $\de$ is a derivation of the Lie algebras
$\OmT$ and $\OmD$\,. Thus $\de$ is compatible with
DGLA structures on $\OmT$ and $\OmD$\,.

Due to the Poincare lemma for the space $\bbRf$
 the complexes $(\OmS,\de)$\,,
$(\OmT,$ $\de)$\,, $(\OmD, \de)$ are acyclic
and their zero cohomologies can be computed
easily, namely

$$
H^0(\OmS,\de)= C^{\infty}(M)\,,
$$
$$
H^0(\OmT,\de)=\FT\,,
$$
and
$$
H^0(\OmD,\de)=\FD\,,
$$
where $\FT$ ($\FD$) denotes the vector space of
all fiberwise po\-ly\-vec\-tor fields (\ref{vect})
(fiberwise po\-ly\-dif\-fe\-ren\-tial operators (\ref{operr}))
with constant coefficients
in $y$'s.

Therefore a natural projection $\si$ from $\OmS$ ($\OmT$, $\OmD$) to
$C^{\infty}(M)$ ($\FT$, $\FD$)
\begin{equation}
\label{si}
\si a= a\Big|_{y^i=dx^i=0}
\end{equation}
gives a morphism of complex $(\OmS,\de)$ ($(\OmT, \de)$\,, $(\OmD, \de)$)
into itself and this morphism is homotopic to the identity map.
One can easily guess the respective homotopy
operator in the from

\begin{equation}
\delta ^{-1}a=y^ki\left( \frac \partial {\partial x^k}\right)
\int\limits_0^1a(x,ty,tdx)\frac{dt}t,  \label{del-1}
\end{equation}
where $a$ is an element of $\OmS$ ($\OmT$\,, $\OmD$)\,,
$i(\partial /\partial x^k)$ stands for the interior
derivative of exterior forms by the vector field
$ \partial /\partial x^k$\,, and $\delta ^{-1}$ is extended to
$C^\infty(M)$ ($\FT$, $\FD$) by zero.

We will use the property of the homotopy
operator $\de^{-1}$ in the following form

\begin{equation}
a=\sigma (a)+\delta \delta ^{-1}a+\delta ^{-1}\delta a\,,
\qquad \forall \quad a\in \Om(M,\cB)\,,
 \label{Hodge}
\end{equation}
where $\cB$ is either $\SM$ or $\T_{poly}$ or $\cD_{poly}$\,.

Thus we have already got a resolution
$(\OmS, \de)$ of the commutative algebra $C^{\infty}(M)$\,.
Now we need to deform this resolution in order to get
$\FT$ and $\FD$ to be identified with $T_{poly}(M)$ and
$D_{poly}(M)$ respectively. In this way we will get resolutions of the
DG Lie algebras $T_{poly}(M)$ and $D_{poly}(M)$\,.

We consider an affine torsion free connection $\n_i$ on $M$
and associate to it the following derivation of $\OmS$

\begin{equation}
\label{nab}
\n= dx^i \frac{\pa}{\pa x^i} + \G
\,:\, \Om^q(M,\SM) \mapsto  \Om^{q+1}(M,\SM)\,,
\end{equation}
where

\begin{equation}
\label{Christ}
\G= -dx^i \G^k_{ij}(x) y^j \frac{\pa}{\pa y^k}\,,
\end{equation}
with $\G^k_{ij}(x)$ being Christoffel symbols of $\n_i$\,.

The derivation $\n$ obviously extends
to derivations of the DG Lie algebras $\OmT$ and $\OmD$

\begin{equation}
\label{nabT}
\n= dx^i \frac{\pa}{\pa x^i} + [\,\G,\, \bul\,]_{SN}
\,:\,  \Om^q (M, \T_{poly}) \mapsto  \Om^{q+1} (M, \T_{poly})\,,
\end{equation}

\begin{equation}
\label{nabD}
\n= dx^i \frac{\pa}{\pa x^i} + [\,\G,\, \bul\,]
\,:\,  \Om^q (M, \cD_{poly}) \mapsto  \Om^{q+1} (M, \cD_{poly})\,.
\end{equation}
It is clear by definition that $\n$ is indeed a derivation of Lie
algebra structures of $\OmT$ and $\OmD$. On the other hand the
multiplication $m\in \G(\cD^1_{poly})$ in $\G(\SM)$ is ``covariantly
constant'' $dm + [\G,m]=0$ and hence the derivation $\n$ commutes
with the differential $\pa$ in $\OmD$\,.

In general derivation (\ref{nab}) is not nilpotent as $\de$. Instead we
have the following expression for $\n^2$

\begin{equation}
\label{nab-sq}
\n^2 a = R \, a
\,:\, \Om^q(M,\SM) \mapsto \Om^{q+2}(M,\SM)\,,
\end{equation}
where
$$
R= -\frac12 dx^i dx^j (R_{ij})^k_l(x) y^l \frac{\pa}{\pa y^k}\,,
$$
and $(R_{ij})^k_l(x)$ is the standard Riemann curvature tensor
of the connection $\n_i$.

Analogously, for $\OmT$ and $\OmD$ we have

\begin{equation}
\label{nab-sq-T}
\n^2 a = [R , a]_{SN}
\,:\,  \Om^q (M, \T_{poly}) \mapsto  \Om^{q+2} (M, \T_{poly})\,,
\end{equation}

\begin{equation}
\label{nab-sq-D}
\n^2 a = [R , a]
\,:\,  \Om^q (M, \cD_{poly}) \mapsto  \Om^{q+2} (M, \cD_{poly})\,.
\end{equation}

Notice that since the connection $\n_i$ is torsion free
derivations $\n$ and $\de$ (anti)commute

\begin{equation}
\label{anticomm}
\de \n + \n \de = 0\,.
\end{equation}

We use the derivation (\ref{nab}) in order to deform the nilpotent
differential $\de$ on $\OmS$, $\OmT$, and $\OmD$\,.

$$
D=\n - \de + A \,:\, \Om^q(M,\SM) \mapsto \Om^{q+1}(M,\SM)\,,
$$
\begin{equation}
\label{DDD}
D=\n - \de + [A\,,\bul\,]_{SN} \,:\,  \Om^q (M, \T_{poly}) \mapsto  \Om^{q+1} (M, \T_{poly})\,,
\end{equation}
$$
D=\n - \de + [A\,,\bul\,] \,:\,  \Om^q (M, \cD_{poly}) \mapsto  \Om^{q+1} (M, \cD_{poly})\,,
$$
where
$$
A=\sum_{p=2}^{\infty}dx^k A^j_{ki_1\ldots i_p}(x) y^{i_1} \ldots
y^{i_p}\frac{\pa}{\pa y^j}
$$
is viewed as an element of $\Om^1(M,\T^0_{poly})$ and an
element of $\Om^1(M,\cD^0_{poly})$\,.

Due to the following theorem one can explicitly construct
a nilpotent differential $D$ in the framework of
ansatz (\ref{DDD})

\begin{teo} Iterating the equation
\begin{equation}
\label{iter_A}
A=\de^{-1} R + \de^{-1}(\n A +\frac12 [A,A])
\end{equation}
in degrees in $y$ one constructs
$A\in \Om^1(M,\T^0_{poly})\subset \Om^1(M,\cD^0_{poly})$
such that $\de^{-1}A=0$ and the derivation $D$ (\ref{DDD})
is nilpotent
$$
D^2=0\,.
$$
\end{teo}
{\bf Proof.} The proof of the theorem is analogous to
the proof of theorem $3.2$ in \cite{Fedosov}\,.

First, we observe that the equation

\begin{equation}
\label{flat}
\de A = R  + \n A +\frac12 [A,A]
\end{equation}
implies that $D^2=0$\,.

Recurrent procedure (\ref{iter_A}) converges to some element
$$A\in \Om^1(M,\T^0_{poly})\subset \Om^1(M,\cD^0_{poly})\,,$$
since the operator $\de^{-1}$ raises the degree in $y$.
An obvious identity $(\de^{-1})^2=0$ implies that for
such $A$
$$\de^{-1}A=0$$
and hence due to homotopy operator property (\ref{Hodge}) we get that
the element $A$ satisfies the
following consequence of equation (\ref{flat})

\begin{equation}
\label{flat1}
\de^{-1} \de A = \de^{-1}R  + \de^{-1}(\n A +\frac12 [A,A])\,.
\end{equation}

We define
$$
C= -\de A + R + \n A +\frac12 [A,A]\,.
$$
Using the Bianchi identities for the Riemann curvature tensor
we get that
$$\de R=0\,,\qquad \n R=0\,.$$
The latter equations imply that $C$ satisfies the following
condition

\begin{equation}
\label{svo}
\n C - \de C + [A, C]=0\,.
\end{equation}
Due to (\ref{flat1}) $\de^{-1}C=0$. Thus
applying $\de^{-1}$ to (\ref{svo}) and using homotopy
operator property (\ref{Hodge}) once again we get that
$$
C= \de^{-1} (\n C + [A,C])\,.
$$
Since the operator $\de^{-1}$ raises the degree in $y$ the
latter equation has a unique zero solution.
Thus the theorem is proved. $\Box$

In what follows we refer to the nilpotent differential $D$ (\ref{DDD})
as Fedosov differential.

Now we are going to prove that the complexes $(\OmS, D)$, $(\OmT,
D)$ and $(\OmD, D)$ are acyclic and their zero cohomologies are
isomorphic to $C^{\infty}(M)$, $\FT$ and $\FD$\,, respectively.

\begin{teo} For a bundle $\cB$ that is either $\SM$
or $\T_{poly}$ or $\cD_{poly}$  we have that
\begin{equation}
\label{h=h0}
H^{\bul}(\Om(M,\cB), D)= H^0 (\Om(M,\cB),D)\,.
\end{equation}
Furthermore,
$$
H^0(\OmS,D)\cong C^{\infty}(M)
$$
as commutative algebras and
$$
H^0(\OmT, D)\cong \FT\,,
$$
$$
H^0(\OmD, D) \cong \FD
$$
as vector spaces.
\label{teo}
\end{teo}
\vspace{0.1cm}
{\bf Proof.} This is a generalization of the proof of
\cite{Fedosov}(Theorem 3.3). Let $a$ be an element of $\Om^q(M,\cB)$ for
$q>0$ and $Da=0$\,. Our purpose is to solve the equation

\begin{equation}
\label{Exact}
a= D b\,.
\end{equation}
We claim that the following recurrent procedure\footnote{For
$\cB=\SM$ one may use $A\,\bul$ instead of $[A,\bul]$\,.}

\begin{equation}
\label{nado}
b= -\de^{-1} a+ \de^{-1}(\n b + [A, b])
\end{equation}
converges to an element $b\in \Om^{q-1}(M, \cB)$ such that
$\de^{-1}b=0$, $\si b=0$, and $Db=a$\,.
All the claims besides the last one are obvious by
construction. Let us prove that $Db=a$\,.

We denote by $h$ the element
$$
h=a-Db \in \Om^q(M,\cB)
$$
and mention that $Dh=0$ or equivalently

\begin{equation}
\label{Closed}
\de h = \n h + [A,h]\,.
\end{equation}

In virtue of equation (\ref{nado})
$$
\de^{-1} h=0\,.
$$
Furthermore, since $q>0$
$$
\si h=0\,,
$$
and hence applying homotopy property (\ref{Hodge}) we get
$$
h= \de^{-1} (\n h + [A,h])\,.
$$
The latter equation has a unique vanishing solution
since $\de^{-1}$ raises the degree in $y$\,.
Thus we have proved (\ref{h=h0}).

We will give the proof only for the iso\-mor\-phism
of the vector spaces
$$
H^0(\OmD, D)\cong \FD
$$
since the analogous statement for $\T_{poly}$
is proved in the same way and the
iso\-mor\-phism of commutative algebras
$H^0(\OmS, D)$ and $C^{\infty}(M)$
is proved in (\cite{Wein}, sect. 6).

As in \cite{Fedosov} we give a constructive proof.
Namely, we will define a bijective map $\tau$ from
$\FD$ to the subspace $\cZ^0(\OmD,D)$ of
$D$-closed forms of degree zero such that for
any $a_0\in \FD$

\begin{equation}
\label{si-tau}
\si (\tau a_0)=\tau a_0\Big|_{y=0}=a_0\,.
\end{equation}
Since $\cZ^0(\OmD, D)=H^0(\OmD,D)$ this would
prove the statement.

For any $a_0\in \FD$ we define an
element $a=\tau a_0\in \Om^{0}(M,\cD_{poly})$
by the following recurrent procedure
\begin{equation}
\label{iter_a}
a=a_0 + \de^{-1}(\n a+ [A,a])\,.
\end{equation}
This procedure converges since $\de^{-1}$ raises the degree in $y$.

First, we prove that $Da=0$ and $\si a= a_0$\,.
While the latter statement is obvious the former one
requires some work.
Let $f=Da$ then $Df=0$, $\si f=0$, and $\de^{-1}f=0$ by (\ref{iter_a}).
Hence due to (\ref{Hodge}) we have
$$
f=\de^{-1}(\n f + [A,f])\,.
$$
This equation has a unique vanishing solution since
$\de^{-1}$ raises the degree in $y$\,.

Thus we have an $\bbR$-linear map $\tau$ from $\FD$ to $\cZ^0(\OmD,
D)$ which is obviously injective. Furthermore,
if $b\in\cZ^0(\OmD, D)$ and $\si b=0$ then due to (\ref{Hodge})
$$
b=\de^{-1}(\n b + [A,b])
$$
and hence $b=0$ since $\de^{-1}$ raises the degree in $y$\,.
Therefore the map $\tau$ is also surjective and
the theorem is proved. $\Box$\\[0.3cm]
{\bf Remark.} The ordinary multiplication $m$ in $\G(\SM)$
viewed as an element of $\G(\cD^1_{poly})$ turns out to be
$D$-closed and
$$
\si m = m \in \FD^1\,,
$$
where $\FD^1$ is a vector space
of fiberwise bidifferential operators on $\SM$ with constant
coefficients in $y$\,. Since $m$ is $D$-closed the Fedosov
differential $D$ (anti)commutes with the differential $\pa$
in $\OmD$. Thus $D$ respects the DGLA structures both on $\OmT$ and
$\OmD$\,.

It turns out that $\FT$ ($\FD$) can be identified with the
vector space $T_{poly}(M)$ ($D_{poly}(M)$) of po\-ly\-vec\-tor fields
(of po\-ly\-dif\-fe\-ren\-tial operators) on $M$\,. Namely, we have
the following

\begin{pred} Given the Fedosov differential (\ref{DDD}) one
can construct an iso\-mor\-phism of vector spaces
\begin{equation}
\label{mu-mu}
\mu\,:\, \FD \mapsto D_{poly}(M)\,,
\end{equation}
whose restriction\footnote{Abusing notations we denote
the iso\-mor\-phism from $\FT$ to $T_{poly}(M)$ by the same letter $\mu$\,.}
to $\FT$ gives an iso\-mor\-phism
from the vector space $\FT$ to the vector space
$T_{poly}(M)$\,. The map $\mu$ preserves degrees of
po\-ly\-dif\-fe\-ren\-tial operators (po\-ly\-vec\-tor fields).
\end{pred}
\vspace{0.1cm}
{\bf Proof.} We give a proof for po\-ly\-dif\-fe\-ren\-tial operators
since a proof for po\-ly\-vec\-tor fields is completely analogous.
First, we restrict ourselves to the
case of degree $0$ po\-ly\-dif\-fe\-ren\-tial operators, that is to
ordinary differential operators.

A construction of the desired
iso\-mor\-phism is based on the observation that for
any function $a_0\in C^{\infty}(M)$
and for any integer $p\ge 0$

\begin{equation}
\label{deriv}
\frac{\pa}{\pa y^{i_1}} \dots \frac{\pa}{\pa y^{i_p}} \tau (a_0)
\Big|_{y=0}= \pa_{x^{i_1}} \ldots \pa_{x^{i_p}} a_0(x)+
{\rm lower~order~ derivatives~ of}~a_0\,.
\end{equation}
Due to this observation an iso\-mor\-phism from $\FD^0$ to $D^0_{poly}$
is defined with the help of the identification between
the space of functions $C^{\infty}(M)$
and the space of $D$-closed sections in $\G(\SM)$\,. Namely,
the following map
$$
\mu\,:\,\FD^0 \,\mapsto\, D^0_{poly}(M)
$$
from the space $\FD^0$ of fiberwise
differential operators on $\SM$ with constant coefficients
in $y$ to the space $D^0_{poly}(M)$ of
differential operators on $M$

\begin{equation}
\label{mu}
\mu(\mP)(a_0)= \si \mP(\tau (a_0))= \mP(\tau (a_0))\Big|_{y=0}\,,
\qquad \mP\in \FD^0\,, \qquad a_0\in C^{\infty}(M)
\end{equation}
gives an iso\-mor\-phism of the respective vector spaces. Then one
can obviously extend the map $\mu$ to the iso\-mor\-phism
from $\FD$ to $D_{poly}(M)$\,. $\Box$

Due to the above proposition $\FT$ ($\FD$) automatically acquires
a structure of DGLA, induced via the map $\mu^{-1}$ from
the vector space $T_{poly}(M)$ ($D_{poly}(M)$).
On the other hand the Fedosov differential $D$
respects the structure of DGLA on $\OmT$
($\OmD$). Thus, the natural question is whether $(\OmT, D)$
($\OmD, D$) is a resolution of $T_{poly}(M)$
($D_{poly}(M)$) as a DGLA or not. The following proposition
gives a positive answer to the question.

\begin{pred}
A DGLA structure induced on cohomologies of the complex $(\OmT,D)$
($(\OmD,D)$) coincides with the DGLA structure induced from $T_{poly}(M)$
($D_{poly}(M)$) via the map $\mu^{-1}$ (\ref{mu-mu})
$$
H^{\bul}(\OmT, D)\cong T_{poly}(M)\,, \qquad
H^{\bul}(\OmD, D)\cong D_{poly}(M)\,.
$$
\end{pred}
\vspace{0.1cm}
{\bf Proof.} As in the previous proof we restrict ourselves to
the case of po\-ly\-dif\-fe\-ren\-tial operators since the proof for
po\-ly\-vec\-tor fields is completely analogous.

First, we prove that the composition  of the maps $\mu$ and $\si$
\begin{equation}
\label{mu-si}
\mu \circ \si \,:\, \cZ^0(\OmD, D)\mapsto D_{poly}(M)
\end{equation}
respects Lie brackets in $\cZ^0(\OmD, D)$ and $D_{poly}(M)$.
Since definitions of the Lie brackets in both Lie algebras
are based on composition of operators (see (\ref{Gerst}))
it suffices to prove that a restriction of the map
(\ref{mu-si}) to the subalgebra $\cZ^0(\Om(M, \cD^0_{poly}), D)$
of $D$-closed elements of $\G(\cD^0_{poly})$ gives a
homomorphism of associative algebras from
$\cZ^0(\Om(M, \cD^0_{poly}), D)$ to the algebra $D^0_{poly}(M)$
of ordinary differential operators on $M$.

To prove this we observe that for any $a_0\in C^{\infty}(M)$ and
for any $D$-closed operator $\mP\in \G(\cD^0_{poly})$  $\mP \tau(a_0)$ is
a D-closed element in $\G(\SM)$. Hence

\begin{equation}
\label{mP-P}
\mP \tau (a_0)= \tau( P a_0)\,,
\end{equation}
where $P=\mu\circ\si(\mP)$\,.

For any two $D$-closed elements $\mP_1$ and $\mP_2$ of
$\G(\cD^0_{poly})$ we denote  $P_1=\mu\circ\si(\mP_1)$ and
$P_2=\mu\circ\si(\mP_2)$\,. Due to (\ref{mP-P}) we have
$$
\tau (P_1P_2a_0)= \mP_1(\tau(P_2 a_0))= \mP_1 \mP_2 \tau(a_0)\,.
$$
Hence
$$
P_1 P_2a_0= \si( \mP_1 \mP_2 \tau(a_0))=
\si(\si (\mP_1\mP_2)\tau a_0) =\mu\circ \si (\mP_1\mP_2) a_0\,.
$$
Thus the map (\ref{mu-si}) respects Lie algebra structures in
$\cZ^0(\OmD, D)$ and $D_{poly}(M)$\,.

Next, we mention that under the map $\mu$ the
multiplication $m\in \FD^1$ in the algebra $\G(\SM)$
turns to the multiplication $m_0\in D^1_{poly}(M)$
in the algebra $C^{\infty}(M)$\,.
$$
\mu (m)=m_0\,.
$$
Due to this observation and the remark made
after the proof of theorem \ref{teo} we have that
for $m\in \cZ^0(\OmD, D)$
$$
\mu\circ\si (m)=m_0\in  D^1_{poly}(M)\,.
$$
Therefore since the differential on $D_{poly}(M)$
is an interior derivation by the element $m_0$ and
the differential on $\cZ^0(\OmD,D)$
is an interior derivation by the element $m$
the map (\ref{mu-si}) is an iso\-mor\-phism
of the DG Lie algebras $\cZ^{0}(\OmD,D)$ and $D_{poly}(M)$\,.
On the other hand we know from theorem \ref{teo} that
$$
H^{\bul}(\OmD,D)\cong \cZ^0(\OmD,D)\,.
$$
Hence the desired statement is proved. $\Box$

\section{Formality theorem for a general ma\-ni\-fold via Fedosov resolutions}
In the previous section we construct Fedosov resolutions
of the algebras of po\-ly\-dif\-fe\-ren\-tial operators and po\-ly\-vec\-tor fields
using the DG Lie algebras $(\OmD,D+\pa,[,])$ and $(\OmT,D,[,]_{SN})$\,.
In terms of strong homotopy Lie algebras this means that we have
two quasi-iso\-mor\-phisms of DG Lie algebras

\begin{equation}
\label{qiT}
U_T\,:\, T_{poly}(M)\leadsto (\OmT,D,[,]_{SN})\,,
\end{equation}

\begin{equation}
\label{qiD}
U_D\,:\, D_{poly}(M)\leadsto (\OmD,D+\pa,[,])\,,
\end{equation}
induced by the homomorphism $\tau\circ \mu^{-1}$\,.
On the other hand we have the fiberwise quasi-iso\-mor\-phism
(\ref{cal-U}) from the DGLA $(\OmT, 0, [,]_{SN})$ to
the DGLA $(\OmD,\pa,[,])$.

In this section we use the quasi-iso\-mor\-phism (\ref{cal-U}) and
the Fedosov resolutions (\ref{qiT}), (\ref{qiD}) in
order to prove that

\begin{teo}[Kontsevich, \cite{Kontsevich}]
\label{ona}
For any smooth ma\-ni\-fold $M$ there exists a quasi-iso\-mor\-phism
$\mmU$ from the DGLA $T_{poly}(M)$ of po\-ly\-vec\-tor fields on $M$ to
the DGLA $D_{poly}(M)$ of po\-ly\-dif\-fe\-ren\-tial operators on $M$.
\end{teo}
~\\
{\bf Proof.} We propose an explicit construction of
the desired quasi-iso\-mor\-phism from the DGLA $T_{poly}(M)$ to
the DGLA $D_{poly}(M)$\,. This construction consists of two steps.

First, we observe that the Fedosov differential (\ref{DDD})
provides us with a Maurer-Cartan element in the DGLA $\OmT$.
By twisting (\ref{cal-U}) with the help of this Maurer-Cartan
element we get the quasi-iso\-mor\-phism

\begin{equation}
\label{qiTD111}
\bU\,:\,(\OmT, D, [,]_{SN})\leadsto (\OmD, D+\pa,[,])\,,
\end{equation}
which readily gives us the quasi-iso\-mor\-phism

\begin{equation}
\label{mU}
\mU\,:\,T_{poly}(M)\leadsto (\OmD, D+\pa,[,])
\end{equation}
due to the presence (\ref{qiT}).
Second, we contract (\ref{mU}) to
the quasi-iso\-mor\-phism

\begin{equation}
\label{mUbar}
\bmU\,:\,T_{poly}(M)\leadsto \cZ^0_{D}(\OmD,\pa,[,])\,,
\end{equation}
which yields the desired quasi-iso\-mor\-phism $\mmU$ from
$T_{poly}(M)$ to $D_{poly}(M)$ since the DG Lie
algebras $D_{poly}(M)$ and $\cZ^0_{D}(\OmD,\pa,[,])$
are isomorphic via the map $\tau\circ \mu^{-1}$\,.

In the remaining part of this section we complete
the proof of theorem \ref{ona} following the
two steps outlined above.

\subsection{Construction of the quasi-iso\-mor\-phism \protect \\
from $(\OmT, D, [,]_{SN})$ to $(\OmD, D+\pa, [,])$}
We present the Fedosov differential (\ref{DDD}) in the
form

$$
D= d + [B, \,\bul\,]_{SN} \,:\,\Om^q(M, \T_{poly}) \mapsto  \Om^{q+1}(M,
\T_{poly})\,,
$$
\begin{equation}
\label{Fed-D}
D= d + [B, \,\bul\,] \,:\, \Om^q(M, \cD_{poly}) \mapsto  \Om^{q+1}(M,
\cD_{poly})\,,
\end{equation}
where
$$d=dx^i \frac{\pa}{\pa x^i}\,,$$
and
\begin{equation}
\label{BBB}
B= -dx^i\frac{\pa}{\pa y^i} -dx^i \G^k_{ij}(x)y^j\frac{\pa}{\pa y^k}
+\sum_{p\ge 2} dx^i A^k_{i j_1\ldots j_p}(x) y^{j_1}\ldots y^{j_p}
\frac{\pa}{\pa y^k}\,.
\end{equation}
Notice that $B$ is only locally viewed as a fiberwise vector
field or a fiberwise differential operator. Namely, if $W$ is
a coordinate disk on $M$ a restriction of $B$ to $W$ gives
an element of $\Om^1(W, \T^0_{poly})$ which can be also viewed
as an element of $\Om^1(W, \cD^0_{poly})$\,. The perhaps surprising thing
is that the transformation law for $B$ upon a change of coordinates
takes a very simple form

\begin{equation}
\label{BBB-tr}
\tilde B\Big|_{\tilde W\cap W}= B\Big|_{\tilde W\cap W}+
dx^i H^k_{ij}(x)y^j\frac{\pa}{\pa y^k}\,.
\end{equation}
Concrete expression for $H^k_{ij}(x)$ is not important. The main
observation we are going to use is that the additional term is
locally a fiberwise po\-ly\-vec\-tor
field, which is linear in $y$'s.

Let us now restrict ourselves to a coordinate disk $W$\,.
On $W$ both the differential $d$ and the element
$B\in \Om^1(W, \T^0_{poly}) \subset \Om^1(W, \cD^0_{poly})$ are
well defined separately. It is easy to see that $d$ commutes with
the fiberwise DGLA structures on $\Om(W,\T_{poly})$ and
$\Om(W,\cD_{poly})$\,.
Moreover, since and the quasi-iso\-mor\-phism $\cU$
(\ref{cal-U}) of the DG Lie algebras $\Om(W,\T_{poly})$, $\Om(W,\cD_{poly})$
is also fiberwise and since $W$ is contractible $\cU$
gives a quasi-iso\-mor\-phism of DG Lie algebras

\begin{equation}
\label{qiTD1}
\cU\,:\,(\Om(W,\T_{poly}), d, [,]_{SN})\leadsto (\Om(W,\cD_{poly}),d+\pa,[,])\,.
\end{equation}

Due to nilpotency of derivations (\ref{Fed-D})
$B\in \Om(W,\T_{poly})\subset \Om(W,\cD_{poly})$
is a Maurer-Cartan element in
the DG Lie algebras $(\Om(W,\T_{poly}), d, [,]_{SN})$ and
$(\Om(W,\cD_{poly}),d+\pa,[,])$\,. Furthermore, using
the terminology of strong homotopy Lie algebras one may say that
the DG Lie algebras $(\Om(W,\T_{poly}), D, [,]_{SN})$
and $(\Om(W,\cD_{poly}), D+\pa, [,])$ are obtained
from the DG Lie algebras $(\Om(W,\T_{poly}), d, [,]_{SN})$
and $(\Om(W,\cD_{poly}), d+\pa, [,])$\,, respectively with
the help of a twisting\footnote{This terminology is borrowed from
\cite{Q} (see App. B $5.3$). However, the twisting by a Maurer-Cartan element
we use here is different from the one in \cite{Q}}
 by the Maurer-Cartan element
$B$\,.
Namely, the nilpotent coderivations $\bQ_T$ and $\bQ_D$
on coassociative cocommutative coalgebras
$C_{\bul}(\Om(W,\T_{poly})[1])$ and $C_{\bul}(\Om(W,\cD_{poly})[1])$
corresponding to the DGLA structures $(D, [,]_{NS})$ and
$(D+\pa, [,])$ on $\Om(W,\T_{poly})$ and $\Om(W,\cD_{poly})$ are
related to the nilpotent coderivations $Q_T$ and $Q_D$\,,
corresponding to the DGLA structures $(d, [,]_{NS})$ and
$(d+\pa, [,])$ as follows

$$\bQ_T(X)=\exp((-B)\wedge) Q_T(\exp(B\wedge) X)\,,$$
\begin{equation}
\label{Q-bQ}
\bQ_D(Y)=\exp((-B)\wedge) Q_D(\exp(B\wedge) Y)\,,
\end{equation}
$$
X\in C_{\bul}(\Om(W,\T_{poly})[1])\,, \qquad Y\in C_{\bul}(\Om(W,\cD_{poly})[1])\,,
$$
where the sums
of the form
$$
\exp(B\wedge) \underbrace{~} =  \underbrace{~} +
B\wedge  \underbrace{~} +\frac{1}{2!} B\wedge B\wedge  \underbrace{~}
+ \ldots
$$
are finite since
$\underbrace{B\wedge B\wedge\ldots \wedge B}_{p}=0$ for
$p > d=dim\, M$\,.

Property $3$ in theorem \ref{1}
implies that the quasi-iso\-mor\-phism $\cU$ (\ref{qiTD1})
maps the Maurer-Cartan element $B$ of
$(\Om(W,\T_{poly}), d, [,]_{SN})$ to $B$\,,
which is viewed as a Maurer-Cartan element of the DGLA
$(\Om(W,\cD_{poly}), d+\pa, [,])$\,.
Therefore a quasi-iso\-mor\-phism from the DGLA
$(\Om(W,\T_{poly}), D, [,]_{SN})$ to the DGLA
$(\Om(W,\cD_{poly}), D+\pa, [,])$ can be obtained by
a twisting the quasi-iso\-mor\-phism $\cU$ with the help
of the Maurer-Cartan element $B$\,. Namely,
the formula

\begin{equation}
\label{cU-bU}
\bU(X)= \exp((-B)\wedge) \cU(\exp(B\wedge)X)\,, \qquad \forall~X\in
C_{\bul}(\Om(W,\T_{poly})[1])
\end{equation}
gives a quasi-iso\-mor\-phism
from the DGLA $(\Om(W,\T_{poly}), D, [,]_{SN})$ to the DGLA
$(\Om(W,\cD_{poly}),$ $D+\pa, [,])$\,. Notice that $B$ which stands
to the right of $\cU$ is viewed
as an element of $\Om(W,\T_{poly})$ and $B$ which
stands to the left of $\cU$ is viewed as an
element of $\Om(W,\cD_{poly})$\,.

Thus we have constructed a quasi-iso\-mor\-phism $\bU$ from
the DGLA $(\Om(W, \T_{poly}), D,[,]_{SN} )$ to the DGLA
$(\Om(W, \cD_{poly}), D + \pa, [,])$
for an arbitrary coordinate disk $W$ on the ma\-ni\-fold $M$\,.
Remarkably, it turns out that the quasi-iso\-mor\-phism $\bU$
does not depend on a choice of local coordinates on $W$
and hence we have the following proposition.

\begin{pred}
Formula (\ref{cU-bU}) defines a quasi-iso\-mor\-phism
from the DGLA $( \Om(M,\T_{poly}),$ $D, [,]_{SN})$ to the DGLA
$(\Om(M,\cD_{poly}),D+\pa , [,])$\,.
\end{pred}
~\\
{\bf Proof.} To prove the assertion we observe that
the $n$-th structure map $\bU_n$ of the quasi-iso\-mor\-phism (\ref{cU-bU})
looks as follows

\begin{equation}
\label{bU-n}
\bU_n(\mv_1, \ldots, \mv_n)= \cU_n(\mv_1, \ldots, \mv_n)+
\sum_{m\ge 1}\frac1{m!}\cU_{n+m}(\underbrace{B,\ldots, B}_{m}, \mv_1, \ldots,
\mv_n)\,,
\end{equation}
$$
\mv_1, \ldots, \mv_n \in  \Om(W,\T_{poly})\,,
$$
where the sum over $m$ is finite since
$\underbrace{B\wedge B\wedge\ldots \wedge B}_{p}=0$ for
$p > d=dim\, M$\,.

Due to property $4$ in theorem \ref{1} and transformation
law (\ref{BBB-tr}) for $B$ the
map $\bU_n$ does not depend on a choice of local coordinates
and therefore $\bU$ is indeed defined as a quasi-iso\-mor\-phism from
the DGLA $(\OmT, D, [,]_{SN})$ to DGLA $(\OmD, D+\pa, [,])$\,. $\Box$

Composing (\ref{qiT}) with $\bU$ we get
a quasi-iso\-mor\-phism

\begin{equation}
\label{mU1}
\mU\,:\,T_{poly}(M)\leadsto (\OmD, D+\pa,[,])\,.
\end{equation}
In the following subsection we use (\ref{mU1})
in order to get a quasi-iso\-mor\-phism from
$T_{poly}(M)$ to $D_{poly}(M)$\,.

\subsection{Contraction of $\mU$ to a quasi-iso\-mor\-phism
from $T_{poly}(M)$ to $D_{poly}(M)$}
In the previous subsection we used a twisting of two DG
Lie algebras and a quasi-iso\-mor\-phism between them. Here we
are going to modify the quasi-iso\-mor\-phism (\ref{mU1})
between the DG Lie algebras $T_{poly}(M)$ and $(\OmD, D+\pa,[,])$
without changing the algebras themselves.
We start with a description of the modification in
a general setting.

Let $(\mh_1,d_1, [\,,\,]_1)$ and
$(\mh_2,d_2, [\,,\,]_2)$ be two DG Lie algebras.
As in section $2$
we associate to  $\mh_1$ and $\mh_2$  coassociative cocommutative
coalgebras $C_{\bul}(\mh_1[1])$ and $C_{\bul}(\mh_2[1])$
with nilpotent coderivations $Q_1$ and $Q_2$\,, induced
by the DGLA structures of $\mh_1$ and
$\mh_2$\,, respectively. Let $F$ be a quasi-iso\-mor\-phism
from $\mh_1$ to $\mh_2$\,.

The statement we are going to use here can be formulated as
\begin{pred} Let
$$
F_m : \wedge^m \mh_1 \mapsto \mh_2[1-m]
$$
be the structure maps of the quasi-iso\-mor\-phism $F$ and $n$ be
any natural number $n\ge 1$\,. Then it is possible
to construct a quasi-iso\-mor\-phism $\tF$ from $\mh_1$ to $\mh_2$
whose structure maps are

\begin{equation}
\label{F-tF}
\cases{
\begin{array}{lr}
\tF_m(\ga_1,\ldots, \ga_m) = F_m(\ga_1, \ldots, \ga_m)\,,
& ~{\rm if }~1\le m < n\,,\\[0.3cm]
\tF_n(\ga_1,\ldots, \ga_n) =
F_n(\ga_1, \ldots, \ga_n)+ d_{2} V_{n}(\ga_1, \ldots, \ga_n) -
 & ~ \\[0.3cm] \displaystyle
\sum_{l=1}^n(-1)^{l+k_1+\ldots +k_{l-1}} 
V_{n}(\ga_1, \ldots, d_1\ga_l, \ldots, \ga_n)   & ~ \\[0.3cm]
\tF_m(\ga_1,\ldots, \ga_m) =
F_m(\ga_1, \ldots, \ga_m)+  W_{m}(\ga_1, \ldots, \ga_n)\,,
& ~ {\rm otherwise,}
\end{array}
}
\end{equation}
where $\ga_l \in \mh_1^{k_l}$\,,
$V_n$ is an arbitrary polylinear map
$$
V_n\,:\,\wedge^n \mh_1 \mapsto \mh_2[-n]
$$
and the maps
$$W_m\,:\,\wedge^m \mh_1 \mapsto \mh_2[1-m]$$
are expressed in terms of $V_n$ and structure maps
$F_m$ via the differentials $d_1$, $d_2$ and
the brackets $[,]_1$ and $[,]_2$\,.
\label{shift}
\end{pred}
~\\
{\bf Proof.} We claim that the desired quasi-iso\-mor\-phism
$\tF$ can be found in the following form

\begin{equation}
\label{tF}
\tF(X) = F(X)+ Q_2 V(X)+ V(Q_1 X)\,,\qquad
\forall~ X\in C_{\bul}(\mh_1[1])\,,
\end{equation}
where $V$ is a linear map from
$C_{\bul}(\mh_1[1])$ to $C_{\bul}(\mh_2[1])$
which satisfies the following relation
$$
\D_2 V(X)= (F\otimes V +  V\otimes F
$$
\begin{equation}
\label{Vrel}
+\frac{1}{2}(V\otimes Q_2 V+Q_2 V \otimes V) +
\frac{1}{2}(V\otimes V Q_1 + V Q_1 \otimes V)) (\D_1 X),
\end{equation}
$$
\forall~ X\in C_{\bul}(\mh_1[1])\,,
$$
$\D_1$ and $\D_2$ denote comultiplications in $C_{\bul}(\mh_1[1])$
and $C_{\bul}(\mh_2[1])$\,, respectively.

The compatibility of the homomorphism (\ref{tF}) with the
nilpotent coderivations $Q_1$ and $Q_2$ follows directly from
the definition
while the compatibility with the comultiplications $\D_1$ and
$\D_2$ follows from equations (\ref{Vrel})\,.

Notice that equations (\ref{Vrel}) simply mean that the linear map
$V$ is defined by structure maps $V_m$ $(m\ge 1)$, which
are arbitrary po\-ly\-li\-near maps

\begin{equation}
\label{V-str}
V_m : \wedge^m \mh_1 \mapsto \mh_2[-m]\,.
\end{equation}
For example, for any $\ga \in \mh_1$
$$V(\ga)=V_1(\ga)\,,$$
and for any pair $\ga_1\in \mh^{k_1}_1$\,, $\ga_2\in \mh^{k_2}_1$
$$ V(\ga_1\wedge\ga_2)=V_2(\ga_1, \ga_2)$$
$$+
(F_1(\ga_1)\wedge V_1(\ga_2)+ \frac{1}{2}V_1(\ga_1) \wedge d_2 V_1(\ga_2)
+\frac{1}{2} V_1(\ga_1)\wedge V_1(d_1 \ga_2)-
(-)^{k_1 k_2}(\ga_1\leftrightarrow\ga_2))\,.
$$
Using (\ref{Vrel}) in this way one can easily express the linear map
$V$ in terms of the structure maps (\ref{V-str}).

Thus the first two structure maps of the shifted
homomorphism $\tF$ takes the from

$$
\tF_1(\ga)=F_1(\ga)+ d_2 V_1(\ga)+ V_1(d_1\ga)\,,
$$
\begin{equation}
\label{tF-str}
\tF_2(\ga_1,\ga_2)= F_2(\ga_1,\ga_2)+ ([F_1(\ga_1), V_1(\ga_2)]_2
+ [V_1(\ga_1),d_2 V_1(\ga_2)]_2
\end{equation}
$$
+[V_1(\ga_1), V_1(d_1 \ga_2)]_2-
(-)^{k_1 k_2}(\ga_1\leftrightarrow\ga_2))-
V_1 ([\ga_1, \ga_2]_1) \,,
$$
where $\ga \in \mh_1$\,, $\ga_1\in \mh^{k_1}_1$ and
$\ga_2\in \mh^{k_2}_1$\,.
The map $\tF_1$ obviously establishes an iso\-mor\-phism of the
spaces of cohomologies $H^{\bul}(\mh_1, d_1)$ and $H^{\bul}(\mh_2,d_2)$
since the map $F_1$ does.

It is clear that the linear map $V$ with a single
non-vanishing structure map
$$
V_n : \wedge^n \mh_1 \mapsto \mh_2[-n]
$$
gives the desired quasi-iso\-mor\-phism $\tF$ and
the proposition follows. $\Box$

Now we are going apply this proposition to the quasi-iso\-mor\-phism
$\mU=\bU\circ U_T$ from the DGLA
$(T_{poly}(M), 0, [,]_{SN})$ to the DGLA
$(\OmD, D+\pa, [,])$ in order to contract it
to a quasi-iso\-mor\-phism $\bmU$ from the
DGLA $(T_{poly}(M), 0, [,]_{SN})$ to the DGLA
$\cZ^0_{D}(\OmD, \pa, [,])$\,.
We formulate a precise statement as

\begin{pred}
One can construct a quasi-iso\-mor\-phism $\bmU$ from the
DGLA $(T_{poly}(M),$ $0, [,]_{SN})$ to the DGLA
$(\OmD, D+\pa, [,])$ such the structure
maps $\bmU_n$ of $\bmU$ take values in
$\cZ^0_{D}(\OmD, \pa, [,])$\,.
\end{pred}
The proof of the proposition is based on the well-known  technique of
spectral sequences. Although the spectral sequence that appears
here is of the simplest type some more care is needed
in using this language since the maps we deal with should
respect the Lie algebra structures as well.
For this reason we give here a detailed
proof.

~\\
{\bf Proof.} We will proceed by induction in $n$\,.\\
{\bf Base of induction.} For $n=1$ we have
\begin{equation}
\label{eee}
(D+\pa) \mU_1(\ga) =0\,, \qquad \forall~\ga\in T_{poly}(M)\,.
\end{equation}
Let
$$
\mU_1(\ga)= \sum_{q=0}^d \mU^q_1(\ga)
$$
be a decomposition of $\mU_1(\ga)$ with respect to the exterior
degree. Then due to (\ref{eee}) (or simply due to
the top degree argument) the component $\mU^d_1(\ga)$ of the maximal
exterior degree $d=dim\,M$ is $D$-closed
$$
D \mU^d_1(\ga)=0\,.
$$
Hence $\mU^d_1(\ga)$ is $D$-exact and there exists a linear map
$V^d_1\,:\,T_{poly}(M)\mapsto \Om^{d-1}(M, \cD_{poly})$
such that
$$
\mU(\ga)+ (D+\pa) V^d_1 (\ga)
$$
has the maximal exterior degree $q_{max}<d$.
Proceeding in this way we can construct
a linear map $V_1\,:\,T_{poly}(M)\mapsto \Om(M, \cD_{poly})$
such that
$$
\mU_1(\ga)+ (D+\pa) V_1 (\ga)
$$
is of exterior degree zero.
Using proposition \ref{shift} we perform the shift
$$
\mU\mapsto \tmU = \mU+ Q_{\OmD}\circ  V + V \circ Q_{T_{poly}(M)}\,,
$$
where the linear map
$V\,:\,C_{\bul}(T_{poly}(M)[1])\mapsto C_{\bul}(\OmD[1])$
has the only non-vanishing structure map $V_1$ constructed
above. This yields a quasi-iso\-mor\-phism $\tmU$ from $T_{poly}(M)$ to
$(\OmD, D+\pa, [,])$ such that for any $\ga\in T_{poly}(M)$
the element $\tmU_1(\ga)$
is of exterior degree zero. Moreover, due to the
equation
$$
(D+\pa) \tmU_1(\ga)=0
$$
$\tmU_1(\ga)$ belongs to $\cZ^0(\OmD, D)$\,.
Thus the base statement of the induction is proved.\\
~\\
{\bf Step of induction.} Let us assume that for all $m<n$ the maps
$\mU_m$ take values in $\cZ^0(\OmD, D)$\,. For the structure map
$\mU_n$ we have

$$
(D+\pa) \mU_n(\ga_1, \ga_2, \ldots, \ga_n) =
$$
\begin{equation}
=\frac12 \sum_{k,l\ge 1,~ k+l=n} \frac1{k!l!}
 \sum_{\ve\in S_n}
\pm [\mU_k (\ga_{\ve_1}, \ldots, \ga_{\ve_k}), \mU_l (\ga_{\ve_{k+1}}, \ldots,
\ga_{\ve_{k+l}})]-
\label{q-iso-mU}
\end{equation}
$$
-\sum_{i\neq j}
\pm \mU_{n-1}([\ga_i,\ga_j]_{SN}, \ga_1, \ldots, \hat{\ga_i}, \ldots, \hat{\ga_j}, \ldots \ga_n),
\qquad \ga_i \in T^{k_i}_{poly}(M)\,.
$$
By the assumption of induction
the right hand side of equation (\ref{q-iso-mU}) is of exterior
degree zero. Hence, by reasoning as above, we
conclude that there exists a polylinear map

$$
V_n\,:\, \wedge^n T_{poly}(M) \mapsto \OmD[-n]\,,
$$
such that for any $\ga_1,\ldots, \ga_n\in T_{poly}(M)$
$$
\mU_n(\ga_1,\ldots, \ga_n) + (D+\pa) V_n(\ga_1,\ldots, \ga_n)
$$
is of exterior degree zero. Therefore, due
to proposition \ref{shift} $\mU$ can be shifted
to a quasi-iso\-mor\-phism $\tmU$\,, whose structure map
$\tmU_n$ takes values in $\Om^0(M,\cD_{poly})$\,.
On the other hand $\tmU_n$ should satisfy the same
cohomological equation (\ref{q-iso-mU}) as $\mU_n$
and hence $\tmU_n$, in fact, takes values
in $\cZ^0(\OmD,D)$\,.
This completes the proof of the proposition and
the proof for theorem \ref{ona}. The desired
quasi-iso\-mor\-phism
$$
\mmU \,:\, T_{poly}(M) \leadsto D_{poly}(M)
$$
is a composition of the quasi-iso\-mor\-phism $\bmU$
and the $DGLA$-iso\-mor\-phism
$$
\mu \circ \si\,:\, \cZ^0_D(\OmD, \pa, [,]) \mapsto D_{poly}(M)\,.~ \Box
$$

As one easily sees, the constructed quasi-iso\-mor\-phism $\mmU$ from the
DGLA $T_{poly}(M)$ to the DGLA $D_{poly}(M)$ depends on a choice of
the affine torsion free connection $\n_i$. In the following
section we will see that symmetries of the connection $\n_i$
determine symmetries of the respective quasi-iso\-mor\-phism $\mmU$\,.

\section{Equivariant formality theorem}
In this section we consider a ma\-ni\-fold $M$ equipped with
a smooth action of a group $G$\,. Given this action on $M$
we can canonically extend it to an action of $G$ on the DGLA $T_{poly}(M)$
of po\-ly\-vec\-tor fields and the DGLA $D_{poly}(M)$ of po\-ly\-dif\-fe\-ren\-tial
operators.  It naturally raises the question as to whether, there
exists a $G$-equivariant formality quasi-iso\-mor\-phism from
the DGLA $T_{poly}(M)$ to the DGLA $D_{poly}(M)$.
Using our construction of the quasi-iso\-mor\-phism of formality
we prove that

\begin{teo}
\label{URA}
If a ma\-ni\-fold $M$ is equipped with a smooth action of
a group $G$ and $M$ admits a $G$-invariant torsion free
connection $\n_i$ then one can construct a $G$-equivariant
quasi-iso\-mor\-phism from the DGLA $T_{poly}(M)$ to
the DGLA $D_{poly}(M)$\,.
\end{teo}
Before giving a proof for the
theorem we mention some interesting corollaries.

First, provided the conditions of theorem \ref{URA} are
satisfied we have a quasi-iso\-mor\-phism between the
respective DG Lie algebras of $G$-invariants

\begin{cor} If a ma\-ni\-fold $M$ is equipped with a
smooth action of a group $G$ and $M$ admits a $G$-invariant
torsion free connection $\n_i$ then one can construct a
quasi-iso\-mor\-phism from the DGLA $(T_{poly}(M))^{G}$ to
the DGLA $(D_{poly}(M))^{G}$\,.
\end{cor}

Second, given a smooth action of a finite or compact group
$G$ one can always construct a $G$-invariant torsion free
connection using the standard averaging procedure.
Hence,

\begin{cor} If a ma\-ni\-fold $M$ is equipped with a
smooth action of a finite or compact group $G$ then one can
construct a $G$-equivariant quasi-iso\-mor\-phism
from the DGLA $T_{poly}(M)$ to the DGLA
$D_{poly}(M)$.
\end{cor}

Since a quasi-iso\-mor\-phism between DG Lie algebras
provides a one-to-one correspondence between the moduli
spaces of Maurer-Cartan elements of the DG Lie algebras
our procedure gives a solution for the deformation quantization
problem of an arbitrary Poisson orbifold. Namely,

\begin{cor} Given a smooth action of a finite group
$G$ on a ma\-ni\-fold $M$ and a $G$-invariant Poisson
structure $\al\in (\wedge^2 T M)^G$ one can always
construct a $G$-invariant star-product $\ast$\,,
corresponding to $\al$\,.
Furthermore, $G$-invariant star-products on $M$ corresponding
to the Poisson bracket $\al$ are classified up to
equivalence by non-trivial $G$-invariant deformations
of $\al$\,.
\end{cor}

Notice that an existence of a star-product on an arbitrary
Poisson orbifold follows in principle from the results
of paper \cite{CFT}. However, in order to get the above
classification of star-products one has to use
equivariant formality theorem \ref{URA}.

Finally, if the group $G$ is neither compact nor finite
the averaging procedure is not applicable. However, one
can still get a $G$-equivariant quasi-iso\-mor\-phism of
formality provided the action of $G$ is ``nice''
enough, namely

\begin{cor} If a ma\-ni\-fold $M$ is equipped with a
free action of a Lie group $G$  then one can
construct a $G$-equivariant quasi-iso\-mor\-phism
from the DGLA $T_{poly}(M)$ to the DGLA
$D_{poly}(M)$.
\end{cor}

\vspace{0.1cm}
{\bf Proof of theorem \ref{URA}.} We are going to prove that
if one uses a $G$-invariant torsion free connection
$\n_i$ in the above construction of the quasi-iso\-mor\-phism
from the algebra $T_{poly}(M)$ of po\-ly\-vec\-tor fields to the
algebra $D_{poly}(M)$ of po\-ly\-dif\-fe\-ren\-tial operators
then the resulting quasi-iso\-mor\-phism is $G$-equivariant.

First, given a smooth action of a group $G$ on $M$ one can
canonically extend the action to the spaces $\OmS$\,,
$\OmT$ and $\OmD$ in such a way that the (super)commutative
product and DGLA structures $(0,[,]_{SN})$ and $(\pa, [,])$
of $\OmS$\,, $\OmT$ and $\OmD$\,, respectively, are
$G$-invariant\,. Second, since $G$ acts on fiber variables
$y$ by linear transformations property $2$ in theorem \ref{1}
implies that the fiberwise
quasi-iso\-mor\-phism (\ref{cal-U}) is $G$-equivariant.

Next, the differential $\de$ (\ref{del}) (\ref{delta}), (\ref{delta1})
the projection $\si$ (\ref{si}) as well as the homotopy operator
$\de^{-1}$ (\ref{del-1}) are obviously $G$-invariant. Moreover,
$G$-invariance of the torsion free connection
$\n_i$ implies $G$-invariance of
the derivation $\n$ (\ref{nab}) (\ref{nabT}), (\ref{nabD})
and $G$-invariance of the
respective Riemann curvature tensor $(R_{ij})^k_l(x)$\,.
Therefore the recurrent procedure (\ref{iter_A}) converges
to a $G$-invariant element $A\in \Om^1(M, \T_{poly})
\subset \Om^1(M, \cD_{poly})$ and hence the respective
Fedosov differential (\ref{DDD}) turns out to be $G$-invariant
as well.

Since the map $\tau$ from $\FD$ ($\FT$) to
$\cZ^0(\OmD,D)$ ($\cZ^0(\OmT,D)$) is defined by
$G$-equivariant recurrent procedure (\ref{iter_a})
it is obviously $G$-invariant. By an analogous line
of arguments we see that the map $\mu$ (\ref{mu-mu}) from $\FD$ ($\FT$) to
$D_{poly}(M)$ ($T_{poly}(M)$) is also $G$-invariant and
therefore quasi-iso\-mor\-phisms (\ref{qiT}) and (\ref{qiD})
are both $G$-equivariant.

Thus it suffices to prove that twisting procedure (\ref{cU-bU})
and the procedure of contraction respect the action of $G$\,.
Since recurrent procedure (\ref{nado}) for solving
cohomological equation (\ref{Exact})
is $G$-equivariant there is nothing to do
with the procedure of contraction. However, some work
is required to prove the $G$-equivariance
of twisting procedure (\ref{cU-bU}).

To show that twisting procedure (\ref{cU-bU}) is $G$-equivariant
we consider the action $\pi(g)$ of $g\in G$ on
the element

\begin{equation}
\label{UBBB}
\mP= \cU_{n+m}(\underbrace{B,\ldots, B}_{m}, \mv_1, \ldots,
\mv_n)\in \OmD\,,
\end{equation}
where $\mv_1, \ldots, \mv_n$
are some elements of $\Om(M,\T_{poly})$\,.
The quasi-iso\-mor\-phism $\bU$ is $G$-equivariant if for any $g\in G$ we
have

\begin{equation}
\label{we-need}
\begin{array}{c}
\pi(g) \cU_{n+m}(\underbrace{B,\ldots, B}_{m}, \mv_1, \ldots,
\mv_n)= ~ \\
 ~ \cU_{n+m}(\underbrace{\pi(g)B,\ldots, \pi(g)B}_{m}, \pi(g)\mv_1, \ldots, \pi(g)\mv_n)\,,
\end{array}
\end{equation}
where $\pi(g)$ acts on all the terms of sum (\ref{BBB})
as on elements of $\Om^1(M,\T^0_{poly})$ besides the
term

\begin{equation}
\label{conn}
\G =-dx^i \G^k_{ij}(x)y^j\frac{\pa}{\pa y^k}\,,
\end{equation}
on which $\pi(g)$ acts by transforming the Christoffel symbols
$\G^k_{ij}(x)$.

On the other hand $G$-equivariance of $\cU$ implies
that

$$
\pi(g) \cU_{n+m}(\underbrace{B,\ldots, B}_{m}, \mv_1, \ldots, \mv_n)=
$$
\begin{equation}
\label{we-have}
 \cU_{n+m}(\underbrace{\pi^{tensor}(g)B,\ldots, \pi^{tensor}(g)B}_{m},
 \pi(g)\mv_1, \ldots, \pi(g)\mv_n)\,,
\end{equation}
where $\pi^{tensor}(g)$ acts on the whole sum (\ref{BBB}) as
on the element of $\Om^1(M,\T^0_{poly})$\,.

But the difference between $\pi^{tensor}(g)B$ and $\pi(g)B$

\begin{equation}
\label{aha}
\pi B - \pi^{tensor}(g) B = H(g)\,.
\end{equation}
is a fiberwise po\-ly\-vec\-tor field linear in the fiber variables $y$'s.
Hence due to property $4$ in theorem \ref{1}
$$
\cU_{n+m}(\underbrace{\pi^{tensor}(g)B,\ldots, \pi^{tensor}(g)B}_{m},
 \pi(g)\mv_1, \ldots, \pi(g)\mv_n)=
$$
$$
\cU_{n+m}(\underbrace{\pi(g)B,\ldots, \pi(g)B}_{m}, \pi(g)\mv_1, \ldots,
 \pi(g)\mv_n)
$$
and the theorem follows. $\Box$

\section{Concluding remarks.}
We conclude the paper by discussing possible generalizations and
applications of the equivariant formality theorem.

Notice that if a Lie group $G$ does not act freely on the ma\-ni\-fold $M$
the existence of a $G$-equivariant formality quasi-iso\-mor\-phism
cannot be guaranteed. However, using the formality  quasi-iso\-mor\-phism
$\mmU$ one can construct a quasi-iso\-mor\-phism $\mmU\,_{\mg}$ from the DGLA of co-chains
$C^{\bul}(\mg, T_{poly}(M))$ to the DGLA of co-chains
$C^{\bul}(\mg, D_{poly}(M))$, where $\mg$ is the Lie algebra
of $G$\,. This quasi-iso\-mor\-phism $\mmU\,_{\mg}$ is
a natural generalization of the quasi-iso\-mor\-phism of
DG Lie algebras $(T_{poly}(M))^G$ and $(D_{poly}(M))^G$ of
invariants. It is interesting to find
relaxed conditions on the action of the Lie group $G$
that would allow us to contract the quasi-iso\-mor\-phism $\mmU\,_{\mg}$
to a quasi-iso\-mor\-phism of the algebras $(T_{poly}(M))^G$ and $(D_{poly}(M))^G$ of
invariants.

An important application of equivariant quasi-iso\-mor\-phisms of
formality is related to quantization \cite{Xu}, \cite{Pasha-N},
\cite{Xu1}, \cite{Takuro}, \cite{Pasha-E}, \cite{Pasha-E1} of classical dynamical $r$-matrices
\cite{Pasha-V}. Despite numerous attempts that have
been undertaken recently in this direction, the
quantization problem of dynamical $r$-matrices remains unsolved
even for the case of triangular dynamical $r$-matrices over
an abelian base. At the moment only certain classes of dynamical
$r$-matrices are known to admit quantization.
Paper \cite{Xu} shows that non-degenerate
triangular dynamical $r$-matrices over
an abelian base can be quantized with the help of the Fedosov
method\footnote{An interesting
modification of the Fedosov method that leads to quantization of
constant triangular $r$-matrices has been proposed in the paper \cite{DILS}}
\cite{Fedosov}.
A class of the so-called completely degenerate
dynamical $r$-matrices has been quantized in
paper \cite{Pasha-N} via the vertex-IRF
transformation \cite{Has}. Interesting examples of
dynamical r-matrices over a non-abelian base found in
works \cite{Pasha-V} and \cite{AM} have been recently quantized
in papers \cite{Pasha-E}, \cite{Pasha-E1}.

In recent paper \cite{Takuro} Kontsevich's formality theorem
has been applied to quantization of triangular dynamical $r$-matrices
over an abelian base. Namely, in \cite{Takuro} it is shown
that if a triangular dynamical $r$-matrix over an abelian base
satisfies a so-called affinization condition then one may
prove a version of equivariant formality theorem
that gives a solution for the quantization problem of
this $r$-matrix.

We would like to mention that quantization of
an arbitrary dynamical $r$-matrix can be
reduced to a version of a problem  of
invariant deformation quantization.
However, in general the action of the respective symmetry
group is not free (in case of non-abelian base
it is not even regular). For this reason we suspect that
there are even examples of triangular dynamical $r$-matrices
over an abelian base which cannot be quantized.

{\bf Acknowledgment.}
This work would not have been possible
without useful discussions with Dmitry Tamarkin to whom I express
my sincere thanks. I also acknowledge Dmitry for his valuable
criticisms concerning the first version of the manuscript.
I would like to thank Pavel Etingof for numerous
stimulating discussions.
I have been benefited from  explanations
by A. Gerasimov, B. Feigin, G. Felder and A. Voronov.
I am grateful to J. Surette for criticisms
concerning my English.
I acknowledge the hospitality of Northwestern University's
Mathematics Department where part of the work
was done. The work is partially supported by the
NSF grant DMS-9988796, the Grant
for Support of Scientific Schools NSh-1999.2003.2,
and the grant INTAS 00-561.

\end{document}